\documentclass[11pt]{article}
\usepackage{amssymb,amsfonts,amsmath}
\newtheorem{theo}{Theorem}

\makeatletter \@addtoreset{equation}{section}
\@addtoreset{theo}{section} \makeatother

\def\qed{\hfill \rule{4pt}{7pt}}
\def\pf{\noindent {\it Proof.} }

\begin{document}


\title{Generalization of matching extensions in graphs (II)
\thanks{This work is supported by RFDP of Higher Education of China and 
Discovery Grant of NSERC of Canada. }}

\author{ Zemin Jin$^1$,  Huifang Yan$^1$ and Qinglin Yu$^1$$^2$\\
\\ {\small $^1$ Center for Combinatorics, LPMC}
\\ {\small Nankai University, Tianjin, PR China}
\\ {\small  $^2$Department of Mathematics and Statistics}
\\ {\small  Thompson Rivers University, Kamloops, BC, Canada}
}

\date{}
\maketitle

\vskip 0.2in

\begin{center}
\begin{minipage}{120mm}

\begin{center}
{\bf Abstract}
\end{center}

{Proposed as a general framework,  Liu and Yu \cite{Liu} ({\it Discrete Math.} 
231 (2001) 311-320)
introduced $(n,k,d)$-graphs to unify the concepts of deficiency of matchings,
$n$-factor-criticality  and $k$-extendability. Let $G$ be a graph and let
$n,k$ and $d$ be non-negative integers such that $n+2k+d\leq
|V(G)|-2$ and $|V(G)|-n-d$ is even. If when deleting any $n$ vertices
from $G$, the remaining subgraph $H$ of $G$ contains a $k$-matching
and each such $k$-matching can be extended to a defect-$d$
matching in $H$, then $G$ is called an
\emph{$(n,k,d)$-graph}. In \cite{Liu}, the recursive relations for distinct
parameters $n, k$ and $d$ were presented and the impact of adding or
deleting an edge also was discussed for the case $d = 0$.  In this paper,
we continue the study begun in \cite{Liu} and obtain new
recursive results for $(n,k,d)$-graphs in the general case $d \geq0$. \\


{\bf AMS(2000) subject classification:} 05C70}
\end{minipage}
\end{center}

\vskip 1in
\baselineskip  0.2in
\newpage

\section{Introduction}

In this paper we consider only finite, undirected and simple
graphs. Denote by $N_G(x)$  set of neighbors of a vertex
$x$ in $G$. If no confusion occurs, we write  $N(x)$ for $N_G(x)$.
Let $G$ be a graph with vertex set $V(G)$ and edge set $E(G)$.
A \emph{matching} $M$ of $G$ is a subset of $E(G)$ such
that any two edges of $M$ have no vertices in common. A matching
of $k$ edges is called a \emph{ $k$-matching}. Let $d$ be a
non-negative integer. A matching is called a \emph{defect-$d$}
matching of $G$ if it covers exactly $|V(G)|-d$ vertices of $G$.
Clearly, a defect-$0$ matching is a perfect matching. A necessary
and sufficient condition for a graph to have a defect-$d$ matching
was given by Berge \cite{Berge}.

\begin{theo} (Berge \cite{Berge})  \label{thmA1}
Let $G$ be a graph and let $d$ be an integer such that $0\leq
d\leq |V(G)|$ and $|V(G)|\equiv d \ (mod \ 2)$. Then $G$ has a
defect-$d$ matching if and only if for any $S\subseteq V(G)$
$$o(G-S)\leq |S|+d.$$
\end{theo}

	For a subset $S$ of $V(G)$, we denote by $G[S]$ the  subgraph of $G$
induced by $S$ and we write $G-S$ for $G[V(G)\setminus S]$. The
number of odd components of $G$ is denoted by $o(G)$. Let $M$ be a
matching of $G$. If there is a matching $M'$ of $G$ such that
$M\subseteq M'$, then we say that $M$ can be extended to $M'$ or
$M'$ is an \emph{extension} of $M$. If each $k$-matching can be extended to a perfect
matching in $G$, then $G$ is called \emph{$k$-extendable}. To avoid triviality, we require
that $|V(G)| \geq 2k+2$ for $k$-extendable graphs.  This family of graphs 
was instroduced by Plummer \cite{Plu80} and studied extensively by Lov\'asz
and Plummer \cite{LP}.  

	A graph $G$ is called \emph{$n$-factor-critical} if after deleting any $n$
vertices the remaining subgraph of $G$ has a perfect matching. 
This concept is introduced by Favaron \cite{Fa96} and Yu \cite{Yu},
independently, which is a generalization of the notions of the well-known
factor-critical graphs and bicritical graphs (the cases of $n = 1$ and $n = 2$).  
Characterizations of $n$-factor-critical graphs, properties of
 $n$-factor-critical graphs and its relationships with other graphic
parameters (e.g., degree sum, toughness, binding number,
connectivity, etc.) have been discussed in \cite{Fa96}, \cite{Fa00} and \cite{Yu}.

	Let $G$ be a graph and let $n,k$ and $d$ be non-negative integers such
that $|V(G)| \geq n+2k+d +2$ and $|V(G)|-n-d$ is even. If when deleting
any $n$ vertices from $G$, the remaining subgraph of $G$ contains
a $k$-matching and each of such $k$-matchings  can be extended to a
defect-$d$ matching in the subgraph, then $G$ is called an
\emph{$(n,k,d)$-graph}. This term was introduced by Liu and Yu \cite{Liu}
as a general framework to unify the
concepts of defect-$d$ matchings, $n$-factor-criticality  and $k$-extendability.
In particular, $(n,0,0)$-graphs are
exactly   $n$-factor-critical graphs and $(0,k,0)$-graphs are just the
same as $k$-extendable graphs. This framework enables the authors to 
prove a series of general results which include many earlier results of matchig theory as 
special cases. In \cite{Liu}, Liu
and Yu provided  the following necessary and sufficient conditions for a graph
to be an $(n,k,d)$-graph.

\begin{theo} \label{thmA2}
A graph $G$ is an $(n,k,d)$-graph if and
only if the following conditions are satisfied.
\begin{enumerate}
\item [(i)] For any $S\subseteq V(G)$ and $|S|\geq n$, then
$$o(G-S)\leq |S|-n+d.$$
\item [(ii)] For any $S\subseteq V(G)$ such that $|S|\geq n+2k$ and
$G[S]$ contains a $k$-matching,
$$o(G-S)\leq |S|-n-2k+d.$$
\end{enumerate}
\end{theo}

	Besides necessary and sufficient conditions, one interesting
problem is to find recursive relations for different parameters
$n, k$ and $d$.  Here, we list some of the relevant results (i.e., Theorems 
\ref{thmA3}-\ref{thmA6}) presented in 
\cite{Liu} for the convenience of the reader.

\begin{theo}  \label{thmA3}

Every $(n,k,d)$-graph $G$ is also an $(n',k',d)$-graph
where $0\leq n'\leq n$, $0\leq k'\leq k$ and $n'\equiv n \ \ (mod\ 2)$.
\end{theo}

In particular, for $d=0$, the following  result was proved.

\begin{theo}  \label{thmA4}
If $G$ is an $(n,k,0)$-graph and $n\geq 1$, $k\geq 2$,
then $G$ is a $(n+2, k-2, 0)$-graph.
\end{theo}

The authors  in \cite{Liu}  also considered   other
recursive properties of $(n,k,d)$-graphs, for instance,
determining the parameters $n^{'}, k^{'}$ and $d^{'}$ such that,
when adding or deleting an edge from an $(n,k,d)$-graph, the
resulting graph is a $(n^{'}, k^{'}, d^{'})$-graph.  The focus in \cite{Liu} 
is mostly on the case of $d=0$ and obtained several
interesting results. For graphs obtained by adding an edge to an
$(n,k,d)$-graph, the following result was shown.

\begin{theo} \label{thmA5}
Let $G$ be an $(n,k,0)$-graph with $n, k\geq 1$. Then
for any edge $e\notin E(G)$, $G\cup e$ is an $(n, k-1, 0)$-graph.
\end{theo}

Moreover, for graphs obtained by deleting an edge from an
$(n,k,d)$-graph, there is the following result.

\begin{theo}  \label{thmA6}
Let $G$ be an $(n,k,0)$-graph, $n\geq 2$ and $k\geq 1$. Then for any edge $e$ of
$G$,
\begin{enumerate}
\item [(i)] $G-e$ is an $(n-2, k, 0)$-graph.
\item [(ii)] $G-e$ is an $(n,k-1,0)$-graph.
\end{enumerate}
\end{theo}

Note that the recursive  results for $d>0$ are not investigated in \cite{Liu}.
In this paper, our main focus is to extend Theorems \ref{thmA4} - \ref{thmA6} to the
case of $d \geq 0$. The results are natural extensions of those in the case of $d = 0$,  but
the proofs are somewhat more involved.
Section 2  is devoted to recursive relations for graphs
obtained by adding an edge to an $(n,k,d)$-graph.  Section 3 presents a 
recursive relation for graphs obtained by adding a vertex. Similar
recursive results for graphs obtained by deleting an edge from an
$(n,k,d)$-graph  are presented in Section 4.

\section{ Recursive relations for adding an edge}

In this section, we consider recursive relations for graphs obtained
by adding an edge to an $(n,k,d)$-graph. First we have the
following result.

\begin{theo} \label{thmB1}
For any $n>d\geq 0$ and $k\geq 1$, if $G$ is an $(n,k,d)$-graph,
then  $G\cup e$ is an $(n,k-1,d)$-graph for any $e \notin E(G)$.
\end{theo}
\pf For $k=1$, since $G$ is an $(n,1,d)$-graph, by Theorem \ref{thmA3},
it is also an $(n,0,d)$-graph.  Hence $G\cup e$ is an
$(n,0,d)$-graph.

So assume that $k\geq 2$.  If $G\cup e$ is not an $(n, k-1,
d)$-graph for some edge $e\notin E(G)$, then there exists an $n$-subset
$S'\subseteq V(G)$ and a $(k-2)$-matching
$M'=\{x_1y_1, x_2y_2, \ldots, x_{k-2}y_{k-2}\}$ such that
the $(k-1)$-matching $e\cup M'$ can not be extended to a
defect-$d$ matching of $G-S'$. Let $e=xy$ and $S'' = V(M')$. By Theorem \ref{thmA1},  
there exists a vertex set $S_1\subseteq G-S'-S''-x-y$ such that
$o(G-S'-S''-x-y-S_1)\geq |S_1|+d+1$.  Since $G$ is  an
$(n,k,d)$-graph, according to Theorem \ref{thmA3}, it is also an $(n,
k-2, d)$-graph. From  Theorem \ref{thmA2} $(ii)$,
$o(G-S'-S''-x-y-S_1)\leq o(G-S'-S''-S_1)+2\leq |S_1|+d+2$. By a
simple parity argument, we have $o(G-S'-S''-x-y-S_1)=|S_1|+d+2$.
Let $S_2=S_1\cup\{x, y\}$. Then, $o(G-S'-S''-S_2)=|S_2|+d$.

\noindent{\it Claim 1.} $S'\cup S_2$ is an independent set in $G$.

Suppose $e_1=uv$ is an edge in $G[S'\cup S_2]$.  Then $uv\cup
M'$ is a $(k-1)$-matching. Let $S=(S'\cup S_2-u-v)\cup (S''\cup
\{u, v\})$ which is of order $|S_2|+n+2(k-1)-2$ and contains a
$(k-1)$-matching. Since $G$ is an $(n,k,d)$-graph,  according to
Theorem \ref{thmA3}, $G$ is also an $(n,k-1,d)$-graph.  Then from 
Theorem \ref{thmA2} $(ii)$ and recall the fact that $|S_2|\geq 2$, we have
$$o(G-S'-S''-S_2)=o(G-S)\leq |S|-n-2(k-1)+d=|S_2|+d-2,$$ a
contradiction.

Let $H = G-S'-S''-S_2$.

\noindent{\it Claim 2.} No even component of $H$ is
connected to $S'\cup S_2$.

Assume that there is an edge, say $e_2= uv$, joining an even
component $C$ of $H$ to $S_2\cup S'$,
where $u\in S'\cup S_2$ and $v\in V(C)$. Then $e_2\cup M'$
is a $(k-1)$-matching. Let $S=(S'\cup S_2-u)\cup (S''\cup\{u,
v\})$ which is of order $n-1+|S_2|+2(k-1)$ and contains a
$(k-1)$-matching. Since $G$ is an $(n,k,d)$-graph,
it is also an $(n, k-1, d)$-graph. Hence Theorem \ref{thmA2} $(ii)$ implies that
$o(G-S)\leq |S|-n-2(k-1)+d=|S_2|-1+d$. However, since the total number of 
odd components increases by at least one upon deleting
$v$ from the even component $C$,  we have that $o(G-S)\geq
o(G-S'-S''-S_2)+1=|S_2|+d+1$,  a contradiction.

\noindent {\it Claim 3.} For every odd component $O$ of $H$, there
do not exist two independent edges $e_3=u_1v_1$ and $e_4=u_2v_2$
joining $O$ to $S'\cup S_2$, where $u_1, u_2\in S'\cup S_2$ and
$v_1, v_2\in V(O)$.

	Suppose, to the contrary, that $e_3$ and $e_4$ are two such 
edges. Then $e_3\cup e_4\cup M'$ is a $k$-matching. Let $S=(S'\cup
S''-u_1- u_2)\cup (S''\cup \{u_1, u_2, v_1, v_2\})$ which is of
order $|S_2|+n-2+2k$ and contains a $k$-matching. Since $G$ is an
$(n,k,d)$-graph, then according to Theorem \ref{thmA2} $(ii)$, we have
$$o(G-S)\leq |S|-n-2k+d=|S_2|+n-2+2k-n-2k+d=|S_2|-2+d.$$
However,  since the total number of odd components  does not
decrease by deleting $v_1$ and $v_2$ from the odd component $O$,
we have $o(G-S)\geq o(G-S'-S''-S_2)=|S_2|+d$, a contradiction.

According to Claim 3, we conclude that for any odd component $O$ of
$H$, if it  is connected to $S_2$ or $S'$ in graph
$G-S''$, then either $|N(V(O))\cap (S'\cup S_2)|=1$ or $|N(S'\cup
S_2)\cap V(O)|=1$.

Since $G$ is an $(n,k,d)$-graph, $G-S''$ is an $(n, 2, d)$-graph
 by Theorem \ref{thmA6} (ii).
Suppose that there are $h$ odd components connected to neither
$S'$ nor $S_2$, and $t$ odd components $C_1, C_2, \ldots, C_t$
with $|N(S'\cup S_2)\cap V(C_i)|=1$,  $1\leq i\leq t$, and
$p=|S_2|+d-h-t$ odd components $D_1, D_2, \ldots, D_{p}$ with
$|N(V(D_i))\cap (S'\cup S_2)|=1$, $1\leq i\leq p $.  Then
$h+t+p=|S_2|+d$. Let $U=\bigcup_{i=1}^{p}{N(V(D_i))\cap (S'\cup
S_2)} =\{u_1, u_2, \ldots, u_q\}$.  We consider the
following three cases:

{\it Case 1}. $n\leq t$. Let $S_3=\bigcup_{i=1}^{n}V(C_i) \cap
N(S'\cup S_2)$. Then $|S_3|=n$. Now we consider the $n$-set
$S_3$ and $(k-2)$-matching $M^{'}$.   From Claim 1, $S'\cup S_2$
is an independent set in $G-S''$. In  $G-S''-S_3$, $S'\cup S_2$
must be matched by vertices of $|S_2|+d-h-n$ odd components from
$C_{n+1}, C_{n+2}, \ldots,C_t, D_1, D_2, \ldots, D_{p}$ and any
maximum matching of $G-S''-S_3$ must miss at least one vertex from
each of $h$ odd components which is connected to neither $S'$ nor
$S''$. Altogether,   a maximum matching of $G-S''-S_3$ will miss
at least
$$
h+|S_2|+n-(|S_2|+d-h-n)=2n+2h-d\geq d+2
$$
vertices (recall that $n>d\geq 0$), which contradicts to the
fact that $G-S''$ is an $(n, 2, d)$-graph.

{\it Case 2.} $t<n\leq q+t$.  Let $S_3=(\bigcup_{i=1}^{t}
V(C_i)\cap N(S'\cup S_2))\bigcup \{u_1, u_2, \ldots, u_{n-t}\}  $.
Now we consider the $n$-set $S_3$ and $(k-2)$-matching
$M^{'}$. Suppose that there are $f$ odd components $D_{i_1},
D_{i_2}, \ldots, D_{i_f}$ among $D_1, D_2, \ldots, D_{p}$ which
are connected to $\{u_1, u_2, \ldots, u_{n-t}\}$ in $G-S''$. It is
obvious that $f\geq n-t$. Note that each vertex of  $(S'\cup
S_2)-S_3$ can only be matched by vertices from $|S_2|+d-h-t-f$ odd
components  $\{D_1, D_2, \ldots, D_{p}\}\setminus \{D_{i_1},
D_{i_2}, \ldots, D_{i_f}\}$ in $G-S''-S_3$. Furthermore,  any
maximum matching of $G-S''-S_3$ must miss at least one vertex from
$D_{i_j}$, $1\leq j\leq f$, and   at least one vertex from each of
$h$ odd components which is connected to neither $S'$ nor $S''$.
Thus any maximum  matching of $G-S''-S_3$ must miss at least
$$
\begin{array}{lll}
f+h+|S_2|+n-(n-t)-(|S_2|+d-h-f-t)&=&2h+2t+2f-d\\
&\geq & 2h+2t+2n-2t-d\\
&\geq & d+2\\
\end{array}
$$
vertices, which implies that $G-S''$ is not an $(n, 2, d)$-graph, a
contradiction again.

{\it Case 3.} $n>q+t$. Let $S_3=(\bigcup_{i=1}^{t} V(C_i)\cap
N(S'\cup S_2))\bigcup U\bigcup S_4  $, where $S_4\subseteq S'\cup
S_2-U$ and $|S_4|=n-q-t$.  Now we consider the $n$-set
$S_3$ and $(k-2)$-matching $M^{'}$.  Note that any maximum
matching of $G-S''-S_3$ must miss at least one vertex from each
of the $h$ odd components connected to neither $S'$ nor $S_2$ and at
least one vertex from $|S_2|+d-h-t$ odd components $D_1, D_2,
\ldots, D_{p}$. Furthermore,  $|S_2|+n-(n-t)$ vertices of $S'\cup
S_2-S_3$ must be missed by any maximum matching of $G-S''-S_3$.
Thus any maximum matching of $G-S''-S_3$  must miss at least
$$
h+|S_2|+d-h-t+|S_2|+n-(n-t)=2|S_2|+d\geq d+4
$$
vertices ($|S_2|\geq 2$),  which implies that $G-S''$ is not an
$(n, 2, d)$-graph, a contradiction again.

This completes the proof.\qed

	Suppose $n, k\geq 1$. Clearly Theorem \ref{thmA5} is a special case of
Theorem \ref{thmB1}. Note that the additional condition $n>d$ in Theorem
\ref{thmB1} is necessary.  For example, consider
a complete bipartite graph $K_{3, d+2}$ with bipartition
$U=\{u_1, u_2, u_3\}$ and $W=\{w_1, w_2, \ldots, w_{d+2}\}$. Let
$H$ be a graph obtained by replacing each $w_i$ by a complete graph $K_{2m+1}$,
$1\leq i\leq d+2$. Obviously, $H$ is a $(1, 2,
d)$-graph, but $H\cup u_1u_2$ is not a $(1, 1, d)$-graph for
$d>0$. An interesting property of the graph $H$ is
that $H$ is a $(1,2,d)$-graph, but not a $(3,0,d)$-graph for
$d>0$. So the conclusion of Theorem \ref{thmA4} does not always hold for $n > d > 0$.

Similarly, under  the additional condition $n>d$, we have the
following result which extends Theorem \ref{thmA4} to the case of $d > 0$.

\begin{theo}\label{thmB2}
For any $n>d\geq 0$ and $k\geq 2$ , if $G$ is an $(n,k,d)$-graph,
then $G$ is also an $(n+2,k-2,d )$-graph.
\end{theo}
\pf Suppose that $G$ is not an $(n+2, k-2, d)$-graph. Then there
exist a vertex set $S'$ of order $n+2$ and ($k-2$)-matching
$M'$ such that $M'$ can not be extended
to a defect-$d$ matching of $G-S'$, i.e., $G-S'-S''$ has no
defect-$d$ matchings.

\noindent {\it Claim.} $S'$ is an independent set in $G$.

If $e=uv$ is an edge in $G[S']$, then $e\cup M'$ can be extended
to a defect-$d$ matching of $G-(S'-u-v)$ since $G$ is an $(n, k-1,
d)$-graph, i.e., $G-S'-V(M')$ has a  defect-$d$ matching, a
contradiction.

Let $u, v$ be two vertices in $S'$ and  $G'=G\cup {uv}$. By
Theorem \ref{thmB1}, $G'$ is an $(n, k-1, d)$-graph. That is,
$uv\cup M'$ can be extended to a defect-$d$ matching $M$ of
$G-(S'-\{u, v\})$. Then $M$ is also a defect-$d$ matching of
$G-S'$ which contains $M'$, a contradiction.

This completes the proof. \qed

\section{Recursive relation for adding a vertex}

Let $G$ be a graph and $x\notin V(G)$.  Denote by $G+x$ the graph
obtained by joining each vertex of $G$ to $x$. Here we consider
the recursive result of adding a vertex to an $(n,k,d)$-graph.

\begin{theo} \label{thmC1}
Let $G$ be an $(n, k, d)$-graph with $k>0$ and $n>d$. Then $G + x$ is
an $(n+1, k-1,d)$-graph for any vertex $x\notin V(G)$.
\end{theo}
\pf Denote $G'=G + x$.  Let $S$ be an $(n+1)$-set of $V(G')$ and
$M'$ a ($k-1$)-matching of $G'-S$. We consider the following
cases:

{\it Case 1.} $x\in S$. Since $G$ is an $(n,k,d)$-graph, it is also
an $(n, k-1, d)$-graph. Let $S'=S-\{x\}$. Then $M'$ can be extended
to a defect-$d$ matching $M$ of $G-S'$ and $M$ is also a
defect-$d$ matching of $G'-S$ which contains the ($k-1$)-matching
$M'$.

{\it Case 2.} $x\in V(M')$. Let $xy$ be an  edge of the
($k-1$)-matching $M'$. If $N(y)\cap S\neq \emptyset$, say $z\in
N(y)\cap S $, then $M''= (M'-xy)\cup yz$ is a ($k-1$)-matching and
$S''=S-\{z\}$ is an $n$-set. Hence $M''$ can be extended to a defect-$d$
matching $M$ of $G-S''$.  It follows that $(M-\{yz\})\cup \{xy\}$ is also a
defect-$d$ matching of $G'-S$ which contains $M'$. If $N(y)\cap S=
\emptyset$, we choose $z$ to be any vertex of $S$. According to
Theorem \ref{thmB1}, $G\cup yz$ is an $(n,k-1,d)$-graph. Since 
$M''=(M'-xy)\cup yz$ be a ($k-1$)-matching and $S''=S-\{z\}$ is an $n$-set, 
$M''$ can be extended to a defect-$d$ matching $M$ of $(G\cup
yz)-S''$.  Then $(M-\{yz\})\cup \{xy\}$ is also a  defect-$d$
matching of $G'-S$ which contains $M'$.

{\it Case 3.} $x\in V(G)-S-V(M')$. Since $G$ is an $(n,k,d)$-graph,
$G$ is  also an $(n,k-1,d)$-graph. Let $y$ be any vertex of $S$ and set $S'=S-y$.
Then $M'$ can be extended to a defect-$d$ matching $M$ of $G-S'$ and 
$d_M(y)=0$ or $d_M(y)=1$. If $d_M(y)=0$, then it is obvious
that $M$ is also a defect-$d$ matching of $G'-S$ which contains
$M'$. If $d_M(y)=1$, let $N_M(y)=z$. Then $(M-yz)\cup xz$ is 
a defect-$d$ matching of $G'-S$. \qed

\section{Recursive relations for deleting an edge}

By presenting an example $H\cong dK_{2m+1}\cup K_2$,
$m\geq 1$, Liu and Yu \cite{Liu} observed that Theorem \ref{thmA6} (i) does
not hold for $d>0$ in general.   Clearly $H$ is a $(2,1,d)$-graph.
But $H-e$ is not a $(0,1,d)$-graph, where $e$ is the edge in the
component $K_2$ of $H$. Furthermore, the graph $H$  implies that Theorem
\ref{thmA6} (ii) does not hold for $d>0$ as well.  Note that the graph $H$
constructed above is not connected. We present a connected example
by modifying $H$ as follows. Let $H'=H + u$.  It is obvious that $H'$
is a $(3, 1,d)$-graph, but $H'-e$ is not a $(1, 1, d)$-graph. Moreover,
$H'$ is a connected counterexample to  Theorem \ref{thmA6} (ii) for $d>0$.

In this section, we provide structural theorems for $G-e$ to be an $(n-2, k, d)$-graph and
an $(n,k-1,d)$-graph, respectively. Also, we discuss the impact of deleting an edge
from bipartite $(n,k,d)$-graphs.

\begin{theo}\label{thmD1}
Let $G$ be an $(n,k,d)$-graph   with  $ n\geq 2$. Then, for an
edge $uv\in E(G)$,  $G-uv$ is not  an  $(n-2, k, d)$-graph if and
only if there exists a vertex subset $S\subseteq V(G)$ with
$|S|=n-2+2k$ such that $G[S]$ contains a $k$-matching and $G-S$
is the union of $d$  odd components, each of which is
factor-critical, and the single edge $uv$.
\end{theo}
\pf $(\Leftarrow)$ The sufficient condition is obvious.

$(\Rightarrow)$  Let $G'=G-uv$.  If $G'$ is not an $(n-2, k,
d)$-graph, then there exists a ($n-2$)-set $S' \subseteq V(G')$
and a $k$-matching $M'$ which can not be extended to a defect-$d$ matching of $G'-S'$.
Let $S'' = V(M')$. Then, by Theorem \ref{thmA1},  there exists a vertex set $S_1\subseteq
V(G')-S'-S''$ such that $o(G'-S'-S''-S_1)\geq |S_1|+d+1$. Then we
have  $\{u, v\}\cap (S'\cup S''\cup S_1)=\emptyset$, for otherwise,
since $G$ is an $(n, k, d)$-graph,  from Theorem \ref{thmA2} $(ii)$, we
have $o(G'-S'-S''-S_1)=o(G-S'-S''-S_1)\leq |S_1|+d$, a
contradiction. Since $G$ is an $(n, k,d)$-graph, we have
$o(G'-S'-S''-S_1)\leq o(G-S'-S''-S_1)+2\leq |S_1|+d+2.$ By a
simple parity argument, we have $o(G'-S'-S''-S_1)= |S_1|+d+2$.
Furthermore,  since $|S_1|+d+2=o(G'-S'-S''-S_1)\leq
o(G-S'-S''-S_1)+2$, we have $o(G-S'-S''-S_1)= |S_1|+d$. Thus $uv$
must be a bridge of an even component of $G-S'-S''-S_1$, which
implies that  $G-S'-S''-S_1$ contains at least one even component.

Let $H = G-S'-S''-S_1$.

\noindent {\it Claim 1.} $H$ has exactly one even component.

Suppose that  $H$ has more than one even  component.
Let $C_1$ and $C_2$ be two such even components of $H$
and $x_1 \in V(C_1)$, $x_2\in V(C_2)$. Since $o(H)=
|S_1|+d$ and, by deleting $x_1$ and $x_2$ from $C_1$ and $C_2$,
the total number of the odd components increases by at least two,
we have $o(H-x_1-x_2)\geq |S_1|+d+2$. However, $G$ is an
$(n, k,d)$-graph, from Theorem \ref{thmA2} $(ii)$, so
$o(G-(S' \cup \{x_1, x_2\})-S''-S_1)=o(H-x_1-x_2)\leq |S_1|+d$, a contradiction.

\noindent {\it Claim 2.}  $|S_1|=0$.

Suppose  $|S_1|\geq 1$. Let $C$ be the even component of
$H$, $x\in S_1$, and $ y\in V(C) $. Since $G$ is an $(n,
k, d)$-graph, from Theorem \ref{thmA2} $(ii)$,  we have
$o(H-y)=o(G-(S'\cup \{x, y\})-S''-(S_1-x))\leq |S_1|+d-1$.
However,  the total number of the odd components increases
when deleting the vertex $y$ from the even component $C$. Since
$o(H)=|S_1|+d$, we have $o(H-y)\geq
|S_1|+d+1$,   a contradiction. Thus $|S_1|=0$.

Let $S=S'\cup S''$. Then $G-S$ is the union of one even component
$C$ which contains edge $uv$ and $d$ odd components $O_1, O_2,
\ldots, O_d$. Since  $o(G'-S'-S''-S_1)= |S_1|+d+2$ and $uv$ is a
bridge of $C$, without loss of generality, we may assume that
$C-uv=O_{d+1}\cup O_{d+2}$. Then $G'-S$ is the union of $d+2$ odd
components $O_1, O_2, \ldots, O_{d+2}$.  Without loss of generality, 
assume $u\in O_{d+1}$ and $v\in O_{d+2}$.

\noindent {\it Claim 3. } $C\cong K_2$ and  each odd component
$O_i$, $1\leq i\leq d$, is factor-critical.

Suppose that  $|V(C)|\geq 4$.   Without loss of generality, assume
that $x$ is a vertex different from $u$ in
$O_{d+1}$. Since $G$ is an $(n,k,d)$-graph, from Theorem
\ref{thmA2} $(ii)$, we have $o(G-(S'\cup \{u, x\})-S'')\leq d$. However, the total number
of the odd components does not decrease by deleting $u$ and $x$
from $O_{d+1}$, which implies that
$o(G-(S'\cup \{u, x\})-S'')=o(G'-(S'\cup \{u, x\})-S'')=d+2$, a contradiction. So
$|V(C)|=2$ and $E(C)=\{uv\}$.

	If $|O_j| = 1$, for all $j$, we are done.  So suppose that for some $j$
($1\leq j\leq d$), $|O_j| \geq 3$ and there exists a vertex $x\in V(O_j)$ such that
$O_j-x$  has no perfect matching. Then any maximum matching of
$G-(S'\cup \{u, x\})-S''$ will miss at least $d+2$ vertices. However, since
$G$ is an $(n, k, d)$-graph,  $G-(S'\cup \{u, x\})-S''$ has a defect-$d$
matching, a contradiction. \qed

From the definition of $(n,k,d)$-graphs, there exists no such
vertex set $S$ mentioned in Theorem \ref{thmD1} for $d = 0$. So Theorem
\ref{thmA6} follows from Theorem \ref{thmD1}.

Though Theorem \ref{thmA6} (i) may not hold for $d>0$ in general, but
there are classes of graphs for which Theorem \ref{thmA6} (i) holds for
$d>0$ without the additional condition $n>d$. We will see that
bipartite graphs are one of such classes.

\begin{theo}\label{thmD2}
Let $G$ be a bipartite  $(n,k,d)$-graph with $ n \geq 2$. Then, for
each edge $e$ of $G$,  $G-e$ is an $(n-2, k, d)$-graph.
\end{theo}
\pf Let $e=uv\in E(G)$. Suppose that  $G-uv$ is not an $(n-2,
k,d)$-graph.  Then, by Theorem \ref{thmD1},  there exists a vertex set
$S\subseteq V(G)$, $|S|=n-2+2k$,  such that $G[S]$ contains a
$k$-matching and $G-S$ is the union of  $d$ factor-critical components and the single edge $e=uv$
since a bipartite graph of order more than $1$ is not factor-critical,
each odd component is a singleton, i.e. $|V(G)|=|S|+d+2=n+2k+d$.
However, from the definition of the  $(n,k,d)$-graph, we have
$n+2k+d\leq |V(G)|-2$, a contradiction.  \qed

	Theorem \ref{thmA6} (ii) does not directly extend to the case $d>0$ in
general. However, sometimes we can characterize the edges which cause the statement in 
Theorem \ref{thmA6} (ii) to fail.

\begin{theo} \label{thmD3}
Let $G$ be an $(n,k,d)$-graph with $k\geq 1$, and  $uv\in E(G)$
such that
$$max\{d_G(u), d_G(v)\}\geq 2k.$$
Then $G-uv$ is  not an $(n, k-1,
d)$-graph if and only if there exists a vertex subset $S\subseteq
V(G)$ with $|S|=n-2+2k$ such that $G[S]$ contains a $(k-1)$-matching
and $G-S$ is the union of $d$ factor-critical odd components and the single edge $uv$.
\end{theo}
\pf  $(\Leftarrow)$ The sufficient condition is obvious.

$(\Rightarrow)$ Let $G'=G-uv$. Suppose that $G'$ is not a
$(n,k-1, d)$-graph. Then there exist a $n$-set $S'\subseteq
V(G)$ and a $(k-1)$-matching $M'$ which  can not be extended to a defect-$d$ matching of
$G'-S'$. Denote $V(M')$ by $S''$. By Theorem \ref{thmA1}, there exists a vertex set $S_1\subseteq
V(G'-S'-S'')$ such that $o(G'-S'-S''-S_1)\geq |S_1|+d+1$. Then we
have  $\{u, v\}\cap (S'\cup S''\cup S_1)=\emptyset$, for otherwise,
since $G$ is an $(n, k, d)$-graph,  from Theorem \ref{thmA2} $(ii)$, we
have $o(G'-S'-S''-S_1)=o(G-S'-S''-S_1)\leq |S_1|+d$, a
contradiction. Moreover, that $G$ is an $(n, k, d)$-graph implies
$o(G'-S'-S''-S_1)\leq o(G-S'-S''-S_1)+2\leq |S_1|+d+2$. By a
simple parity argument, we conclude $o(G'-S'-S''-S_1)=|S_1|+d+2$ and
$o(G-S'-S''-S_1)=|S_1|+d$. Thus $uv$ must be a bridge of an even
component $C$ of $G-S'-S''-S_1$, which implies that
$G-S'-S''-S_1$ contains at least one even component.

\noindent {\it Claim 1.} $((N_G(u)\cup N_G(v))\cap
(V(G)-S'-S''))-\{u, v\}=\emptyset$.

Suppose that  $ux$ is an edge in $G-S'-S''-v$. Since $G$ is an
$(n,k,d)$-graph, $ux\cup M'$ is a $k$-matching of $G-S'$ which can
be extended to a defect-$d$ matching $M$ of $G-S'$. Then $M$ is a
defect-$d$ matching which contains $M'$ but  not $uv$, a
contradiction.

Claim  1 implies that $C$ is a complete graph consisting of the
single edge $uv$.

\noindent {\it Claim 2.} $S_1=\emptyset$.

Without loss of generality, assume that  $d_G(u)\geq 2k$ (i.e.,
$d_G(u)>|S''|+|\{v\}|$). Thus $N(u)\cap S'\neq \emptyset$ or
$N(u)\cap S_1 \neq \emptyset$. Consider the case of $N(u)\cap
S'\neq \emptyset$. Let $x \in N(u)\cap S'$ and $y\in S_1\neq
\emptyset$.  Since $G$ is an $(n,k,d)$-graph, the $k$-matching
$M'\cup ux$ can be extended to a defect-$d$ matching of $G-(S'\cup
y - x)$. Thus $o(G-(S'\cup y-x)-(S''\cup ux) -(S_1-y))\leq
|S_1|-1+d$. On the other hand, since $o(G-S'-S''-S_1)=|S_1|+d$ and
$C$ is a single edge, $G-(S'\cup y-x)-(S''\cup ux)-(S_1-y)$ has
$|S_1|+d+1$ odd components, a contradiction. For the case of
$N(u)\cap S_1 \neq \emptyset$, we obtain a similar contradiction.

\noindent {\it Claim 3.} $C$ is the only even component of
$G-S'-S''$.

The arguments are similar to that of Claim 2. Suppose that there
is another even component $C'$ in $G-S'-S''$. Let $y\in V(C')$.
Then there exists an edge $ux \in E(C, S')$ so that the
$k$-matching $M' \cup ux$ can be extended to a defect-$d$ matching
of $G-(S'\cup y-x)$ which implies that $o(G-(S'\cup y-x)-(S''\cup
ux)-S_1)\leq |S_1|+d$. However,  since $o(G-S'-S''-S_1)=|S_1|+d$
and the number of odd components increases upon deleting $y$ from
$C'$,  $G-(S'\cup y-x)-(S''\cup ux)-S_1$ has at least $|S_1|+d+2$
odd components, a contradiction.

\noindent {\it Claim 4.} Each odd component of $G-S'-S''$ is
factor-critical.

Suppose that $O$ is an odd component of $G-S'-S''$ which is not
factor-critical. Hence there exists a vertex $y\in V(O)$ such that
$O-y$ has no perfect matching. Since $G$ is an $(n, k,d)$-graph,
$G-S''$ is an $(n, 1, d)$-graph. Thus, for any $x \in N_G(u) \cap
S'$, $ux$ can be extended to a defect-$d$ matching of $G-(S'\cup
y-x)-S''$, which is impossible since such a matching will miss at least $d+2$
vertices.

Let $S=S'\cup S''$. From the claims above, $G-S$ is  the union of
$d$ factor-critical odd components and a single edge $uv$. \qed

Finally,  we present an example to show that the condition 
$max\{d_G(u), d_G(v)\}\geq 2k$ in Theorem \ref{thmD3} is necessary. Let $G$
be the graph with vertices $x_1, x_2, x_3, x_4, x_5$ and the edges
$x_1x_2, x_2 x_3, x_3x_4, x_4x_5, x_5x_1, x_2x_4,
x_3x_5$. Taking $n$ disjoint copies of $G$ and an edge $e = uv$,  
join the vertices $u$ and $v$ to $x_3$ and $x_4$ in each copy of $G$.
Denote the resulting graph by $H$. Then $max\{d_H(u),
d_H(v)\}=2n+1<2(n+1)$.  One can verify that $H$ is an
$(1,n+1,n+1)$-graph and $H-uv$ is not an $(1,n,n+1)$-graph.
However, for any vertex subset $S\subseteq V(H)$ with $|S|=2n+1$
such that $H[S]$ contains a $n$-matching, $H-S$ is not the union
of $n+1$ factor-critical odd components and a single edge $uv$.

	This article is merely the first of series of investigations of a general 
framework to unify the various extendabilities and factor-criticalities.  So far we have discussed 
the characterization of $(n, k, d)$-graphs and the recursive relations only. The important aspects of 
$(n, k, d)$-graphs, such as decomposition procedure, Gallai-type structural theorems and algorithms 
for finding $(n, k, d)$-graphs, have not been explored yet.  More research on this subject will follow.

\vskip 1mm \vspace{0.3cm}

\noindent \title{\Large\bf Acknowledgments}
\maketitle

The authors are indebted to the anonymous referees for providing a detailed 
comments and suggestions. 

\small

\end{document}